\documentclass{article}

\newtheorem{thm}{Theorem}[section]

\newtheorem{ex}[thm]{Example}

\newcommand{\be}{\begin{equation}}
\newcommand{\ee}{\end{equation}}
\newcommand{\ben}{\begin{enumerate}}
\newcommand{\een}{\end{enumerate}}
\newcommand{\pa}{{\partial}}
\newcommand{\R}{{\rm R}}
\newcommand{\e}{{\epsilon}}
\newcommand{\g}{{\bf g}}
\newcommand{\pxi}{{\pa \over \pa x^i}}

\newcommand{\qed}{\hspace*{\fill}Q.E.D.}  
\title{\Large On the Flag Curvature of  Finsler Metrics \\
of Scalar Curvature }
\author{Xinyue Chen\footnote{supported by the National Natural Science Foundation of China (10171117)}, Xiaohuan Mo\footnote{supported by the National Natural Science Foundation of China (10171002)} and Zhongmin Shen}

\date{March 15, 2003}

\begin{document}

\maketitle

\begin{abstract} The flag curvature of a Finsler metric is  called a Riemannian quantity  because it is an extension of sectional curvature in Riemannian geometry. 
In Finsler geometry, there are several  non-Riemannian quantities such as the (mean) Cartan torsion, the (mean) Landsberg curvature and the S-curvature, which all vanish for Riemannian metrics. It is important to understand the geometric meanings of these quantities.
In this paper, we study  Finsler metrics of scalar curvature (i.e., the flag curvature is a scalar function on the slit tangent bundle) and   partially determine the flag curvature  when certain non-Riemannian quantities are isotropic. 
 Using the obtained formula for the flag curvature,  we classify locally projectively flat Randers metrics  with isotropic S-curvature. 
\end{abstract}

\bigskip

\section{Introduction}

Finsler metrics arise naturally from many areas of mathematics as well as natural science. For example, the navigation problem  in a Riemannian  space  gives rise to a lots of interesting Finsler metrics with special geometric properties  \cite{Sh5} \cite{Sh6} \cite{Zer}. 
In Finsler geometry, we  study not only the shape of a space, but also the ``color'' of the space on an infinitesimal scale. The Riemannian quantity (such as the flag curvature) describes the shape of a space, while non-Riemannian quantities describes the ``color'' of the space.

For a Finsler manifold $(M, F)$, the  flag curvature ${\bf K}={\bf K}(P, y)$ is a function of  tangent planes $P={\rm span}\{ y, v\}\subset T_xM$ and directions $y\in P\setminus\{0\}$. This quantity tells us how curved the space is at a point. If $F$ is  Riemannian, ${\bf K}={\bf K}(P)$ is independent of $y\in P\setminus\{0\}$, ${\bf K}$ being called the sectional curvature in 
Riemannian geometry.  
A Finsler metric $F$ is  said to be {\it of scalar curvature} if the flag curvature 
${\bf K}= {\bf K}(x, y)$ is a scalar function on the slit tangent bundle $TM\setminus\{0\}$. 
Clearly, a Riemannian metric is of scalar curvature if and only if ${\bf K}={\bf K}(x)$ is a scalar function on $M$ (which is a constant in dimension $n>2$ by the Schur lemma).
There are lots of non-Riemannian Finsler metrics of scalar curvature.
  One of the important problems in Finsler geometry is to study and characterize Finsler metrics of scalar curvature. This problem has not been solved yet, even for Finsler metrics of constant flag curvature.

 According to E. Cartan's local classification theorem,  any  Riemannian metric $\alpha$  of constant sectional curvature $\mu$ is locally isometric to the following standard metric 
$\alpha_{\mu}$  on the unit ball ${\rm B^n}\subset \R^n$ or the whole $\R^n$ for $\mu=-1, 0, +1$: 
\begin{eqnarray}
\alpha_{-1} (x, y) & = & { \sqrt{|y|^2-(|x|^2|y|^2-\langle x, y \rangle^2)}\over 1-|x|^2},  \ \  \ \ \ y\in T_x {\rm B}^n\cong \R^n, \label{alpha-1}\\
\alpha_0(x, y) & = & |y|,\ \ \ \ \ \  \ \ \ \ \ \ \ \ \ \ \ \ \ \  \ \ \ \ \ \ \ \  \ \ \ \ \ \ \ \  y\in T_x\R^n \cong \R^n,  \label{alpha-0}\\
 \alpha_{+1}(x,y)& = & {\sqrt{|y|^2+(|x|^2|y|^2-\langle x, y \rangle^2)}\over 1+|x|^2}\ \ \ \ \ \  y\in T_x\R^n \cong \R^n.\label{alpha+1}
\end{eqnarray}

The simplest non-Riemannian Finsler metrics are those in the form $F=\alpha+\beta$,  where $\alpha $ is a Riemannian metric and $\beta$ is a $1$-form.  They are called {\it Randers metrics}.  Bao-Robles prove that 
if a Randers metric $F= \alpha+\beta$ has isotropic flag curvature ${\bf K}={\bf K}(x)$,  then there is a constant $c$ such that the covariant derivatives of $\beta$ with respect to $\alpha$ satisfy a system of PDEs \cite{BaRo}
(i.e.,  equation (\ref{e_00}) below with $c(x)=c$). Recently, Bao-Robles-Shen have classified 
Randers metrics of constant curvature via the navigation problem  in Riemannian manifolds \cite{BaRoSh}. 

A Finsler metric is said to be {\it locally projectively flat} if at any point there is a local coordinate system  in which the geodesics are straight lines as point sets. Why are we interested in these type of Finsler metrics?
Riemannian metrics of constant curvature are locally projectively flat. The converse is true too according to Beltrami's theorem.  Projectively flat 
Finsler metrics on a convex domain in $\R^n$ are regular solutions to Hilbert's Fourth Problem \cite{Hi}. It is known that every
locally projectively flat Finsler metric is  of scalar curvature.    Locally projectively flat Finsler metrics with constant flag curvature have been  solved at satisfactory level \cite{Be1} \cite{Be2}
\cite{Br1}-\cite{Br3}, \cite{Fk1}-\cite{Fk3},  \cite{Sh3},  \cite{Sh4}.

In Finsler geometry,  there are several important non-Riemannian quantities: the distortion  $\tau$,  the mean Cartan torsion ${\bf I}$, the S-curvature ${\bf S}$ and the mean Landsberg curvature ${\bf J}$, etc. They all vanish  for Riemannian metrics, hence they are said to be {\it non-Riemannian}. See  Section \ref{sectionP}
for more details about their definitions and geometric meanings.

All known Randers metrics $F=\alpha+\beta$  of scalar curvature (in dimension $n>2$) satisfy 
${\bf S} = (n+1)c(x) F$ or ${\bf J} + c(x) F {\bf I}=0$, where $c(x)$ is a scalar curvature (see Theorem \ref{ChenShen_thm1.1} below for projectively flat examples). In order to classify Finsler metrics of scalar curvature, we first investigate those with isotropic S-curvature.

\begin{thm}\label{thmKS}
Let $(M, F)$ be an $n$-dimensional    Finsler manifold of scalar curvature with flag curvature ${\bf K}(x, y)$. Suppose that the S-curvature is isotropic, 
\be
{\bf S} = (n+1) c(x)  F(x, y),\label{Sc}
\ee
where $c(x)$ is a scalar function on $M$. Then there is a scalar function $\sigma(x)$ on $M$ such that
\be
{\bf K} =3 {c_{x^m}(x)y^m\over F(x, y)} + \sigma(x). \label{KFc}
\ee
In particular,  $c(x)=c$ is a constant if and only if 
${\bf K} = {\bf K}(x)$ is a scalar function on $M$. 
\end{thm}

\bigskip

In (\ref{KFc})  and thereafter, the subscript $x^m$ in $c_{x^m}$  indicates partial differentiation with respect to $x^m$. 
In Theorem \ref{thmKS}, we partially determine the flag curvature when the S-curvature is isotropic. This  is a generalization of a theorem in \cite{Mo2} where the second author  shows that the flag curvature is isotropic,  ${\bf K}={\bf K}(x)$ if  (\ref{Sc}) holds for $c(x)=constant$.
In this case, ${\bf K}= constant$ when $n \geq 3$ by the Schur theorem \cite{BCS}.

\begin{thm}\label{thmKtau} Let $(M, F)$ be an $n$-dimensional  Finsler manifold of scalar curvature. Suppose that  ${\bf J}/{\bf I}$ is isotropic, \be
{\bf J}+ c(x) F{\bf I}=0, \label{JcF}
\ee
where $c =c(x)$ is a $C^{\infty}$ scalar function on $M$. Then the flag curvature 
${\bf K}={\bf K}(x, y)$ and the distortion $\tau=\tau(x,y)$ satisfy
\be
{n+1 \over 3} {\bf K}_{y^k} + \Big ( {\bf K} + c(x)^2 - { c_{x^m}(x) y^m \over F(x,y)} \Big ) \tau_{y^k} =0.\label{KKc}
\ee
\ben
\item[(a)] If $c(x)= c$ is a constant, then  there is a scalar function  $\rho(x)$ on $M$ such that 
\be
{\bf K} = - c^2 + \rho(x) e^{- {3 \tau(x, y)\over n+1} }, \ \ \ \ \ \ y\in T_xM\setminus\{0\}. \label{eq6}
\ee
\item[(b)] Suppose that $F$ is non-Riemannian on any open subset of $M$. If  ${\bf K}={\bf K}(x)$ is a scalar function on $M$, then  $c(x)=c$ is a constant, in which case ${\bf K}= -c^2\leq 0$. \een
\end{thm}

\bigskip
From Theorem \ref{thmKtau}(a), it seems that (\ref{JcF}) is weaker than (\ref{Sc}) when $c(x)=constant$. But for a non-Riemannian Randers metric $F=\alpha+\beta$, (\ref{JcF}) implies (\ref{Sc}) \cite{ChSh}.

Finsler metrics with ${\bf J}=0$ are said to be {\it weakly Landsbergian}. Berwald metrics are weakly Landsbergian.
According to Theorem \ref{thmKtau}(a), for any weak Landsberg metric  of scalar curvature, 
${\bf K} = \rho(x) \exp ( - {3 \tau(x, y)\over n+1} )$. 
 Here is an open problem: Is there any  weak Landsberg metric of scalar curvature, which is non-Berwaldian? 
Further, one would like to  know whether or not there are any non-Berwaldian Finsler metric of scalar curvature satisfying 
${\bf J}+ c F {\bf I}=0$? The only known examples are the 
 family of (generalized) Funk metrics on the unit ball ${\rm B}^n \subset \R^n$  satisfying 
$ {\bf J}+{1\over 2} F {\bf I}=0$ and ${\bf K}= -{1\over 4}$  (see Theorem \ref{ChenShen_thm1.1} (A2) below).

Theorem  \ref{thmKtau}(a) in the case when $c=0$ is essentially proved in Matsumoto's book. See Proposition
26.2 in \cite{Ma}. Matsumoto assumes that $F$ is a Landsberg metric, but what  he actually  needs in his proof is that ${\bf J}=0$. Since the notion of distortion has not been introduced in \cite{Ma} by that time, Matsumoto's proposition is stated in a local coordinate system.

\bigskip

As we have mentioned early, locally projectively flat Finsler metrics are of scalar curvature. 
For a Randers metric $F=\alpha+\beta$, it is locally projectively flat and only if 
$\alpha$ is locally projectively flat (equivalently, of constant sectional curvature by Beltrami's theorem) and $\beta$ is closed \cite{BaMa} \cite{MoSh}, in which case, (\ref{Sc}) holds if and only if (\ref{JcF}) holds \cite{ChSh}. A natural problem  is to determine $\beta$ such that  $F=\alpha+\beta$ is locally projectively flat with ${\bf S}= (n+1)c(x)F$. 
Using  Theorem \ref{thmKS}, we find an explicit expression for $\beta$.
 
\begin{thm}\label{ChenShen_thm1.1}
Let $F = \alpha+\beta$ be a locally projectively flat Randers metric on an $n$-dimensional manifold $M$ and $\mu$ denote the constant sectional curvature of $\alpha$. Suppose that the S-curvature is 
isotropic,
${\bf S}= (n+1) c(x) F$. Then $F$ can be classified as follows. 
\ben
\item[(A)] If  $\mu + 4 c(x)^2 \equiv 0$, then $c(x)=constant$ and  the flag curvature ${\bf K} = - c^2$.
\ben 
\item[(A1)] if $c=0$, then $F$ is locally Minkowskian with flag curvature ${\bf K}=0$;
\item[(A2)] if $c \not=0$, then after a normalization, $F$ is locally isometric to the following    
Randers metric on the unit ball ${\rm B}^n \subset \R^n$,
\be
F(x,y)= {\sqrt{ |y|^2 - ( |x|^2|y|^2-\langle x, y\rangle^2)} \pm \langle x, y \rangle 
\over 1- |x|^2 } \pm {\langle a, y \rangle \over 1+\langle a, x \rangle},\label{F_a}
\ee
where $a\in \R^n$ with $|a|<1$, and the flag curvature of $F$ is negative constant, 
${\bf K} = - {1\over 4}$. 
\een
\item[(B)]If  $ \mu + 4 c(x)^2 \not=0$, then $F$  is  given by
\be
F(x, y) = \alpha(x, y)  - { 2 c_{x^k} (x) y^k \over \mu + 4 c(x)^2 }\label{eq5}
\ee
and the flag curvature ${\bf K}$ of $F$ is given by
\[
{\bf K} = 3 \Big \{ {c_{x^k}(x)y^k\over F(x,y)} + c(x)^2 \Big \} + \mu
= {3\over 4}\Big \{ \mu + 4 c(x)^2 \Big \} {F(x, -y)\over F(x, y)} + {\mu\over 4}.
\]
\ben
\item[(B1)] when $\mu=-1$, $\alpha=\alpha_{-1}$ can be expressed in the form  (\ref{alpha-1}) on ${\rm B}^n$. In this case,  
\[
c(x) = { \lambda + \langle a, x \rangle \over 2 \sqrt{ (\lambda + \langle a, x \rangle )^2 \pm (1- |x|^2 ) } },
\]
where $\lambda\in \R$ and $a\in \R^n$ with $ |a|^2 < \lambda^2 \pm 1 $.
\item[(B2)] when $\mu=0$, $\alpha=\alpha_0$ can be expressed in the form (\ref{alpha-0}) on $\R^n$. In this case, 
\[
c(x) = {\pm 1\over 2 \sqrt{ k + 2\langle a, x\rangle + |x|^2 } },
\]
where $k >0$ and $a\in \R^n$ with $|a|^2 < k$.
\item[(B3)] when $\mu = 1$, $\alpha=\alpha_{+1}$ can be expressed in the form (\ref{alpha+1}) on $\R^n$.
 In this case,
\[
c(x) = {\e+ \langle a, x \rangle  \over 2 \sqrt{ 1 + |x|^2 - (\e + \langle a, x \rangle )^2 } },
\] 
where $\e\in \R$ and $a\in \R^n$ with $ |\e|^2 +|a|^2 < 1$.
\een
\een 
\end{thm}

\bigskip
 In a forthcoming paper, we are going to determine the local structures of projectively flat Finsler metrics with isotropic S-curvature.

Theorem \ref{ChenShen_thm1.1} is a local classification theorem. If we assume that the manifold is compact without boundary, then 
the scalar function $c(x)$ takes much more special values.

\begin{thm}\label{ChenShen_thm1.2}
Let $F=\alpha+\beta$ be a  locally projectively flat Randers metric on an $n$-dimensional compact manifold $M$ without boundary. Let $\mu$ denote the constant sectional curvature of $\alpha$. Suppose that ${\bf S}= (n+1) c(x) F$.
\ben
\item[(a)] If 
$\mu =-1$, then 
$F =\alpha$ is Riemannian.
\item[(b)] If $\mu=0$, then $F$ is locally Minkowskian.
\item[(c)] If $\mu = 1$, then $c(x) ={1\over 2} f(x)/\sqrt{1-f(x)^2}$ and 
\be
F(x,y)= \alpha(x,y) - {f_{x^k}(x)y^k \over \sqrt{1-  f(x)^2 } },\label{Fmu=0}
\ee
where $f(x)$ is an eigenfunction of the standard Laplacian corresponding to the eigenvalue $\lambda=n$ with 
$\max_{x\in M} |f|(x) <1$.  
Moreover, the flag curvature and the S-curvature of $F$ are  given by 
\be
{\bf K}(x,y) = {3\over 4 (1-f(x)^2)} 
{F(x, -y)\over F(x, y) } + {1\over 4}. \label{KFF}
\ee
\[
{\bf S}(x,y) = { (n+1) f(x) \over 2 \sqrt{ 1- f(x)^2} }  F (x, y).
\]
\een
\end{thm}

In the case when $\mu=1$, $(M, \alpha)$ is isometric to the lense space ${\rm S}^n/\Gamma$.
Let  $F=\alpha+\beta$ be the Randers metric on ${\rm S}^n$ defined  in (\ref{Fmu=0}) using some eigenfunction $f$ on ${\rm S}^n$.  
It can shown that 
$\delta:= \sqrt{|\nabla f|^2_{\alpha}(x) + f(x)^2} <1$ is a constant. 
By (\ref{KFF}), we obtain the following bounds on the flag curvature of $F$.
\[
{2-\delta \over 2(1+\delta)}  \leq {\bf K} \leq  { 2 +\delta \over 2 (1-\delta)}.
\]
Assume that $(M, \alpha) ={\rm S}^n$ is the standard unit sphere. 
Since $F$ is pointwise projectively equivalent to $\alpha$,  
the geodesics of $F$ are great circles. One can easily see that the $F$-length of any great circle is equal to  $2\pi$.

\section{Preliminaries}\label{sectionP}

In this section, we are going to 
give a brief description on several geometric quantities in Finsler geometry.

Let $F$ be a Finsler metric on an $n$-dimensional manifold $M$. 
The geodesics of $F$ are characterized by the following equations
\[
\ddot{c}^i(t) + 2 G^i \Big (c(t), \dot{c}(t) \Big ) =0,
\]
where $G^i=G^i(x, y)$ are given by
\[ G^i = {1\over 4} g^{il}\Big \{ [F^2]_{x^ky^l} y^k - [F^2]_{x^l} \Big \}.\]
where $g_{ij}(x,y) = {1\over 2} [F^2]_{y^iy^j}(x,y)$ and $(g^{ij}(x,y)):=(g_{ij}(x,y))^{-1}$.
When $F$ is Riemannian, i.e., $g_{ij}(x,y)= g_{ij}(x)$ depend only on $x\in M$,
$G^i(x,y) = {1\over 2} \Gamma^i_{jk}(x)y^jy^k$ are quadratic in $y= y^i\pxi|_x$. 
There are many non-Riemannian Finsler metrics with this property. Such Finsler metrics are called {\it Berwald metrics}.  

Let 
\[ \tau (x, y):=\ln \Big [ { \sqrt{\det \Big ( g_{ij}(x,y) \Big )}\over {\rm Vol} ({\rm B}^n(1)) }
 \cdot {\rm Vol} \Big \{ (y^i) \in \R^n \Big | \ F\Big (y^i \pxi|_x \Big ) < 1 \Big \} \Big ] .\]
$\tau=\tau(x,y)$ is a scalar function on $TM\setminus\{0\}$, which is called  the  {\it distortion} \cite{Sh2}.

Let 
\[ I_i(x,y): = {\pa \tau \over \pa y^j}(x,y)=
 {1\over 2} g^{jk}(x,y) {\pa g_{jk}\over \pa y^i}(x,y).\]
The quantity 
${\bf I}_y := I_i(x,y)dx^i$  is called the {\it mean Cartan torsion}. 
According to Deicke's theorem,  $F_x$ is Euclidean at $x\in M$ if and only if ${\bf I}_y=0$, or equivalently, $\tau=\tau(x)$  at $x\in M$ \cite{BCS} \cite{De}. 

Let 
\[ {\bf S}(x, y):= { d \over dt} \Big [ \tau \Big ( \sigma(t), \dot{\sigma}(t) \Big ) \Big ]_{t=0},
\]
where $\sigma(t)$ is the geodesic with $\sigma(0)=x$ and $\dot{\sigma}(0)=y$. 
${\bf S}$ is called the {\it S-curvature} \cite{Sh1}\cite{Sh2}. 
There are lots of Randers metrics of constant flag curvature satisfying ${\bf S}=(n+1) c F$ for some constant $c$ \cite{BaSh} \cite{Sh2}  \cite{Sh5} \cite{Sh6}.
${\bf S}$ said to be  {\it isotropic} if there is a scalar functions 
$c(x)$ on $M$ such that 
\[ {\bf S}(x, y) = (n+1) c(x) F(x, y).\]

The horizontal covariant derivatives of ${\bf I}$ along geodesics  give rise to  the mean Landsberg curvature ${\bf J}_y: = J_i (x,y) dx^i$, 
where $J_i=J_i(x,y)$ are given by
\[ J_i: = y^m {\pa I_i \over \pa x^m} -I_m {\pa G^m\over  \pa y^i} - 2 G^m {\pa I_i \over \pa y^m}.\]
A Finsler metric $F$ is said to be {\it weakly Landsbergian} if ${\bf J}=0$. ${\bf J}/{\bf I}$ is regarded as the relative rate of change of ${\bf I}$ along geodesics. The generalized Funk metrics on the unit ball ${\rm B}^n \subset \R^n$ satisfy ${\bf J}+ c F {\bf I}=0$ for some constant $c \not=0$ \cite{ChSh} \cite{Sh3}. 
 ${\bf J}/{\bf I}$  is said to be {\it isotropic} if there is a scalar function $c(x)$ on $M$ such that
\[ {\bf J}+ c(x) F {\bf I}=0.\]

\bigskip

The 
Riemann curvature ${\bf K}_y= K^i_{\ k}  dx^k \otimes \pxi|_x : T_xM \to T_xM$ is a family of linear maps on tangent spaces,
defined by
\be
K^i_{\ k} = 2 {\pa G^i\over \pa x^k}-y^j{\pa^2 G^i\over \pa x^j\pa y^k}
+2G^j {\pa^2 G^i \over \pa y^j \pa y^k} - {\pa G^i \over \pa y^j}{\pa G^j \over \pa y^k}.  \label{Riemann}
\ee
For a flag $P={\rm span}\{y, u\} \subset T_xM$ with flagpole $y$, the {\it flag curvature} ${\bf K}={\bf K}(P, y)$ is defined by
\[
{\bf K}(P, y):= {\g_y (u, {\bf K}_y(u))
\over \g_y(y, y) \g_y(u,u)
-\g_y(y, u)^2 },
\]
where $\g_y = g_{ij}(x,y)dx^i\otimes dx^j$. 
 When $F$ is Riemannian, ${\bf K}={\bf  K}(P)$ is independent of $y\in P$, which  is just
the sectional curvature of $P$ in Riemannian geometry. 
We say that  a Finsler metric $F$ is  {\it of scalar curvature} if 
for any $y\in T_xM$, 
the flag curvature ${\bf K}= {\bf K}(x, y)$ is a scalar function on the slit tangent bundle $TM\setminus\{0\}$.
If ${\bf K}=constant $, then $F$ is said to be of {\it constant flag  curvature}.

It is easy to see that 
$F$ is locally projectively flat if and only if at any point there is a standard local coordinate system $(x^i, y^i)$ in $TM$ such that $G^i(x, y)= P(x, y)y^i$. In this case,
it  follows from (\ref{Riemann}) that $K^i_{\ k} = \Xi \delta^i_k + \tau_k \; y^i $. Thus 
$F$ is of scalar curvature.

\bigskip
Let $M$ be an $n$-dimensional manifold. 
Let $\pi^*TM$ denote the pull-back tangent bundle by 
$\pi: TM\setminus\{0\}\to M$ and $(x, y, v)$ denote the elements of $\pi^*TM$, where $y\in T_xM \setminus\{0\}$ and $v\in T_xM$. 
Let $\pi^*T^*M$ denote the horizontal cotangent bundle of $TM\setminus\{0\}$,
consisting of $\pi^*\theta$, where $\theta\in T^*M$. There is a natural duality between $\pi^*TM$ and $\pi^*T^*M$. 
Let $\{ {\bf e}_i:=  (x, y, \pxi)\}$ be a natural local frame for $\pi^*TM$. Then   $\{ \omega^i:= \pi^*dx^i\}$ is the  dual local coframe for $\pi^*T^*M$.
$\pi^*TM$ has a canonical section,  $Y:=(x, y, y)=y^i {\bf e}_i$,  where $y=y^i\pxi$. 

Given a Finsler metric $F$ on $M$. It  defines the Riemannian metric tensor ${\bf g}= g_{ij}\omega^i\otimes \omega^j$ and the  Cartan torsion ${\bf C}=C_{ijk} \omega^i\otimes \omega^j \otimes \omega^k$ on $\pi^*TM$, where $g_{ij}= {1\over 2} [F^2]_{y^iy^j}$ and $C_{ijk}= {1\over 4} [F^2]_{y^iy^jy^k}$. 
 The Chern connection is a linear connection on $\pi^*TM$, which  are characterized by 
\[ d\omega^i = \omega^j \wedge \omega_j^{\ i},\]
\[dg_{ij} = g_{ik} \omega_j^{\ k} + g_{kj} \omega_i^{\ k} + 2 C_{ijk} \{   dy^k + y^j \omega_j^{\ k}\}.
\]
 See \cite{BCS}\cite{Ch}. 
Let
\[\omega^{n+k} := dy^k + y^j \omega_j^{\ k}.\]
We obtain a local coframe $\{ \omega^i, \omega^{n+i}\}$ for $T^*(TM\setminus\{0\})$. 
Let 
\[ \Omega^i:= d\omega^{n+i}- \omega^{n+j} \wedge \omega_j^{\ i}.\]
We can express $\Omega^i$ in the following form
\[
\Omega^i = {1\over 2} K^i_{\ kl}\omega^k \wedge \omega^l - L^i_{\ kl}\omega^{k} \wedge \omega^{n+l},
\]
where $K^i_{\ kl} + K^i_{\ lk}=0$. 
Let 
\[
K^i_{\ k}:= K^i_{\ kl}y^l.
\]
We obtain the Riemann curvature 
${\bf K}= K^i_{\ k}\omega^k \otimes \omega^l$ 
and 
${\bf L}=L^i_{\ kl} \omega^k \otimes \omega^l \otimes {\bf e}_i$. 
In a standard local coordinate system $(x^i, y^i)$, $K^i_{\ k}$ are given by (\ref{Riemann}). 
The Riemann curvature is introduced by Riemann in 1854 for Riemannian metrics and extended to Finsler metrics by L. Berwald in 1926 \cite{Be1}\cite{Be2}.

With the Chern connection, we define covariant derivatives of quantities on $TM\setminus\{0\}$  in the usual way. For example, for a scalar function $f$, we define $f_{|i}$ and $ f_{\cdot i}$ by
\[ df = f_{|i} \omega^i + f_{\cdot i} \omega^{n+i},\]
for the mean Cartan torsion ${\bf I}=I_i \omega^i$, define $I_{i|j}$ and $I_{i\cdot j}$ by
\[ d I_i - I_k \omega_i^{\ k} = I_{i|j}\omega^j + I_{i\cdot j} \omega^{n+j}.\]
For a tensor ${\bf T}= T_{i \cdots k } \omega^i \otimes \cdots \otimes \omega^k$, 
\[  T_{i \cdots k \; \cdot m} = {\pa T_{i\cdots k}\over \pa y^m}.\]
Without much difficulty, one can show that
\[
K^i_{\ kl} = {1\over 3} \Big \{ K^i_{\ k\cdot l} - K^i_{\ l \cdot k} \Big \}
\]
and
\be
I_i = \tau_{\cdot i}, \ \ \ \ \ {\bf S}:= \tau_{|m}y^m, \ \ \ \ \
 J_i = I_{i|m}y^m. \label{ILJ}
\ee
Moreover,  $L_{ijk}:=g_{im}L^m_{\ kl} = C_{ijk|m}y^m$. See \cite{Sh2}. The following equations are proved in \cite{Mo1} \cite{MoSh}.
\begin{eqnarray}
L_{ijk|m}y^m + C_{ijm}K^m_{\ k} & = &  - {1\over 3}g_{im}K^m_{\ k\cdot j}
- {1\over 3} g_{jm} K^m_{\ k\cdot i}\nonumber\\
& &  - {1\over 6} g_{im}K^m_{\ j\cdot k}
- {1\over 6} g_{jm}K^m_{\ i\cdot k}. \label{Moeq1}
\end{eqnarray} 
Contracting (\ref{Moeq1}) with $g^{ij}$
gives
\be
J_{k|m}y^m + I_mK^m_{\ k}  = -  {1\over 3}\Big \{ 2 K^m_{\ k\cdot m} +  K^m_{\ m\cdot k}\Big \} . \label{Moeq2}
\ee

As a scalar function on $TM\setminus\{0\}$, the distortion satisfies the following Ricci identities
\be
\tau_{|k|l} = \tau_{|l|k} + \tau_{\cdot m} K^m_{\ kl},\label{tauK}
\ee
\be
\tau_{|k\cdot l} = \tau_{\cdot l | k} -\tau_{\cdot m} L^m_{\ kl}. \label{TauRicci}
\ee
Contracting (\ref{TauRicci}) with $y^k$ yields
\be
 {\bf S}_{\cdot l} = \tau_{|l} +  J_l. \label{tauSI}
\ee
It follows from (\ref{tauSI}) that
\be
{\bf S}_{\cdot k |l} = \tau_{|k|l} + J_{k|l}. \label{StauJ}
\ee
Using (\ref{ILJ}), (\ref{tauK}),   and (\ref{StauJ}), we obtain
\begin{eqnarray*}
{\bf S}_{\cdot k|l}y^l - {\bf S}_{| k} & = & \Big ( {\bf S}_{\cdot k |l} - {\bf S}_{\cdot l | k} \Big )y^l \\
& = & \Big ( \tau_{|k|l} - \tau_{|l|k} \Big ) y^l
+ \Big ( J_{k|l} - J_{l|k} \Big ) y^l \\
& = &  \tau_{\cdot m} K^m_{\ \ kl}y^l  - I_m K^m_{\ \ k} -  {1\over 3}\Big \{ 2 K^m_{\ k\cdot m} +  K^m_{\ m\cdot k}\Big \}\\
& = & -  {1\over 3}\Big \{ 2 K^m_{\ k\cdot m} +  K^m_{\ m\cdot k}\Big \}.
\end{eqnarray*}
We obtain 
\be
{\bf S}_{\cdot k|m}y^m - {\bf S}_{| k} 
= - {1\over 3} \Big \{ 2 K^m_{\ k\cdot m} + K^m_{\ m\cdot k} \Big \}.\label{SSKKMo}
\ee
Equation (\ref{SSKKMo}) is established in \cite{Mo2}.

Now we assume that $F$ is of scalar curvature with flag curvature ${\bf K}= {\bf K}(x, y)$. This is equivalent to the following identity:
\be
K^i_{\ k} = {\bf K} F^2 \; h^i_k, \label{Kikiso1}
\ee
where $h^i_k := g^{ij} h_{jk}$ and $h_{jk}:= g_{jk} -F^{-2} g_{js}y^s g_{kt}y^t$. 
By  (\ref{Moeq1}),  (\ref{Moeq2}) and (\ref{Kikiso1}), we obtain 
\[
 L_{ijk|m}y^m =  - {1\over 3}F^2 \Big \{ {\bf K}_{\cdot i} h_{jk}   + {\bf K}_{\cdot j} h_{ik}  + {\bf K}_{\cdot k} h_{ij} + 3 {\bf K} C_{ijk}\Big \}
\]
and
\be
J_{k|m}y^m = - {1\over 3}F^2\Big \{ (n+1) {\bf K}_{\cdot k} + 3{\bf K} I_k \Big \}.\label{AZeq2}
\ee

\section{Proofs of Theorems \ref{thmKS} and \ref{thmKtau}}
In this section, we are going to prove the first two theorems.

\bigskip
\noindent
{\it Proof of Theorem \ref{thmKS}}: 
Plugging  (\ref{Kikiso1})  into (\ref{SSKKMo}), we obtain 
\be
{\bf S}_{\cdot k|l}y^l - {\bf S}_{| k} 
= - {n+1\over 3} {\bf K}_{\cdot k} F^2.\label{SSKF}
\ee
Plugging (\ref{Sc}) into (\ref{SSKF}) yields
\be
c_{|l}(x)y^l F_{\cdot k} - c_{|k}(x) F = - {1\over 3} {\bf K}_{\cdot k} F^2.\label{cFcF}
\ee
It follows from (\ref{cFcF}) that 
\[
\Big [   {1\over 3} {\bf K} - { c_{|m}(x) y^m \over F(x, y) } \Big ]_{y^k}=0.
\]
Thus 
\[
\sigma:= {\bf K} - {3 c_{|m}y^m \over F} 
\]
is a scalar function on $M$. 
This proves the theorem.
\qed

\bigskip

\bigskip
\noindent
{\it Proof of Theorem \ref{thmKtau}}: By assumption, $J_k = - c F I_k$ and $J_k = I_{k|m}y^m$ we obtain 
\[
 J_{k|m}y^m = - c_{|m}y^m F I_k - cF I_{k|m}y^m 
= - c_{|m}y^m F I_k + c^2 F^2 I_k.
\]
It follows from (\ref{AZeq2}) that 
\be {n+1 \over 3} {\bf K}_{\cdot k} + \Big ( {\bf K} + c^2 - { c_{x^m} y^m \over F} \Big ) I_k =0. \label{KKcI}
\ee
By (\ref{ILJ}),  $I_k = \tau_{\cdot k}$. We obtain (\ref{KKc}).

(a) Suppose that $c_{x^m}(x)= 0$ at some point $x\in M$. Then equation (\ref{KKc}) simplifies to
\[
{n+1 \over 3} {\bf K}_{y^k} + \Big ( {\bf K} + c^2  \Big ) \tau_{y^k} =0.
\]
This implies that 
\[ \Big [ \Big ({\bf K}+c^2\Big )^{n+1\over 3} e^{\tau} \Big ]_{y^k}
=\Big ( {\bf K}+ c^2 \Big )^{n-2\over 3} e^{\tau} \Big \{ {n+1\over 3}
{\bf K}_{y^k} + \Big ({\bf K}+c^2 \Big ) \tau_{y^k} \Big \} =0.
\]
Thus the function $({\bf K}+c^2)^{n+1\over 3} e^{\tau} $ is independent of $y\in T_xM$. There is a number $\rho(x)$ such that 
\be
 {\bf K} = - c(x)^2 + \rho(x) e^{-{3\tau(x,y)\over n+1} }. \label{Kcrho}
\ee
 When  $c(x)= c$ is a constant, we obtain (\ref{eq6}) from (\ref{Kcrho}).
Note that $\rho(x)$ is not necessarily a constant. 

(b) Suppose that ${\bf K}={\bf K}(x)$ is a scalar function on $M$. Then 
(\ref{KKc}) simplifies to
\be \Big ( {\bf K} + c^2 - { c_{x^m} y^m \over F} \Big ) \tau_{y^k} =0. \label{KKcI*}
\ee
We claim that $c(x) = c$ is a constant. If this is false, then there is an open subset ${\cal U}$ such that  $dc(x) \not=0$ for any $x\in {\cal U}$.
Clearly, at any $x\in {\cal U}$, 
 ${\bf K}(x) \not= -c(x)^2 + c_{x^m}(x)y^m/F(x,y)$ for almost all $y\in T_{x}M$. By (\ref{KKcI*}),
$\tau_{\cdot k}= I_k =0$. Thus $F$ is Riemannian on ${\cal U}$ by Deicke's theorem (cf. \cite{De} \cite{BCS}). 
 This contradicts our assumption in the theorem. 
This proves the claim. By (\ref{Kcrho})
and  (\ref{KKcI*}), we obtain
\be
\rho(x)\; \tau_{y^k} =0. \label{KKcI**}
\ee
We claim that $\rho(x)\equiv 0$. If this is false, then there is an open 
subset ${\cal U}$ such that $\rho(x)\not=0$ for any $x\in {\cal U}$. 
By (\ref{KKcI**}), we obtain that $\tau_{y^k} = I_k =0$ on ${\cal U}$. Thus $F$ is Riemannian on ${\cal U}$. This again contradicts the assumption 
in the theorem. Therefore $\rho(x)\equiv 0$. We conclude that
${\bf K} = -c^2$ by (\ref{Kcrho}).  
\qed

\bigskip
According to \cite{ChSh}, for any Randers metric $F=\alpha+\beta$,  (\ref{JcF}) holds if and only if (\ref{Sc}) holds and $\beta$ is closed.
For a general Finsler metric, (\ref{JcF}) does not imply (\ref{Sc}). 
Now we combine two conditions (\ref{Sc}) and (\ref{JcF}) and prove the following

\begin{thm}
\label{thmSJ} Let $(M, F)$ be an $n$-dimensional  Finsler manifold of scalar curvature.  Suppose that the S-curvature and    the  mean Landsberg curvature satisfy
\be
{\bf S} = (n+1) c(x) F, \ \ \ \ \ {\bf J} + c(x) F {\bf I}=0, \label{SJ}
\ee
where $c =c(x)$ is a scalar function on $M$.
Then the flag curvature is given by
\be
 {\bf K} = 3 {c_{x^m}(x) y^m \over F(x,y)} + \sigma(x) = - { 3c(x)^2+\sigma(x)\over 2} + \nu(x) e^{ - 2\tau(x,y)\over n+1},\label{Kmu}
\ee
where $\sigma(x)$ and $\nu(x)$ are  scalar functions on $M$.
\ben 
\item[(a)] Suppose that  $F$ is not Riemannian on any open subset in $M$. 
If $c(x)=c$ is a constant, then 
${\bf  K} = - c^2, $ $\sigma(x) = -c^2$ and $\nu(x)=0$.
\item[(b)] If $ c(x) \not=constant$, then
the distortion is given by
\be
\tau = \ln \Big \{ \frac{2\nu(x)F(x,y)}{6 c_{x^m}(x)y^m + 3[\sigma(x) + c(x)^2 ]F(x,y)}
\Big \}^{2\over n+1}.     \label{tauF}
\ee
\een
\end{thm} 
{\it Proof}:
By the above argument, ${\bf K}$ is given by (\ref{KFc})
and it satisfies (\ref{KKc}). It follows from (\ref{KFc}) that
\[ {c_{x^m}(x) y^m \over F(x, y)} = {1\over 3} \Big ( {\bf K} - \sigma(x)\Big ).\]
Plugging it into (\ref{KKc}) yields
\[
{n+1\over 3} {\bf K}_{y^k} + \Big ( {2\over 3} {\bf K} +  c(x)^2 +{1\over 3}\sigma(x) \Big ) \tau_{y^k} =0.
\]
We obtain 
\[ \Big [ \Big ( 2 {\bf K} + 3c(x)^2 +\sigma(x) \Big )^{n+1\over 2} e^{\tau} \Big ]_{y^k}=0.
\]
Thus there is a scalar function $\nu(x)$ on $M$ such that 
\be
 {\bf K} = - { 3c(x)^2 +\sigma(x) \over 2} + \nu (x) e^{ - {2\tau(x, y) \over n+1} }.\label{Kv}
\ee
Comparing (\ref{Kv}) with (\ref{KFc}), we obtain 
\be
{c_{x^m}(x)y^m
\over F(x, y)} 
= - { c(x)^2+\sigma(x)\over 2} + {\nu(x)\over 3} e^{- {2\tau(x, y)\over n+1} }. \label{ctau}
\ee

(a) Suppose that  $c(x)=c$ is a constant. We claim that $\nu(x)=0$. If it false, then
${\cal U}:=\{ x\in M, \ \nu(x) \not=0\} \not=\emptyset$. From (\ref{ctau}), one can see that  $\tau=\tau(x)$ is a scalar function on  ${\cal U}$, hence $F$ is Riemannian on ${\cal U}$ by Deicke's theorem \cite{De} \cite{BCS}. This contradicts the assumption in (a). 

Now  (\ref{ctau}) is reduced to that $\sigma(x)  = c(x)^2$ and (\ref{Kv})  is reduced to that  ${\bf K}=-c^2$.

(b) If $c(x)\not=contant$, then $\nu(x)\not=0$ by (\ref{ctau}). In this case, we can solve 
(\ref{ctau}) for $\tau$ and obtain (\ref{tauF}).
\qed

\bigskip

It follows from Theorem \ref{thmSJ} that if a Finsler metric of scalar curvature satisfies that  
${\bf S} = (n+1) c F$ and ${\bf J}+ c F {\bf I}=0$ for some constant $c$, then the flag curvature is given by  ${\bf K}=- c^2$. 
One would like to know  whether or not there are non-Riemannian, non-locally Minkowskian  Finsler metrics with these properties. If a Randers metric has these properties, then it is, up to a scaling, locally isometric to the generalized Funk metric on the unit ball ${\rm B}^n \subset \R^n$ \cite{ChSh}. 
 In dimension two, any Finsler metric with ${\bf S}=0, {\bf J}=0$ and ${\bf K}=0$ is locally Minkowskian.

\bigskip

\begin{ex}\label{exJ+cI=0}{\rm   For an arbitrary number $\e$ with $0 < \e \leq 1$, 
define
\begin{eqnarray*}
\alpha: & = & {\sqrt{ (1-\e^2)(xu+yv)^2 +\e (u^2+v^2) (1+ \e (x^2+y^2) )}\over 1 + \e (x^2 +y^2) }\\
\beta : & = & {\sqrt{1-\e^2}(xu+yv)\over 1+ \e (x^2+y^2) }.
\end{eqnarray*}
We have
\[ \|\beta\|_{\alpha} =\sqrt{1-\e^2}\sqrt{ x^2+y^2 \over \e + x^2+y^2 } < 1.\]
Thus $F:=\alpha+\beta$ is a Randers metric on $\R^2$. In \cite{ChSh}, we have verified  that 
\[{\bf S} = 3 c F,  \ \ \ \  {\bf J}_y + c F \; {\bf I}_y =0\]
where
\[ c = { \sqrt{1-\e^2}\over 2 ( \e + x^2+y^2 ) },\]
and obtained a formula for the Gauss curvature
\begin{eqnarray*}
 {\bf K}
& = &  { - 3 \sqrt{1-\e^2}\; (xu+yv)/(1+\e (x^2+y^2)) 
\over  \sqrt{ (1-\e^2)(xu+yv)^2 +\e (u^2+v^2) (1+ \e (x^2+y^2) )}+ \sqrt{1-\e^2}(xu+yv)      }\\
&& + { 7 (1-\e^2) + 8\e (\e + x^2 +y^2 )\over 4 ( \e + x^2+y^2 )^2 }.
\end{eqnarray*}

\bigskip
Here we are going to compute $\sigma$ and $\nu$ in Theorem \ref{thmSJ}. 
By a direct computation we can express the function $\sigma := {\bf K}-{3(c_xu+c_yv)\over F} $ in (\ref{Kmu})  by
\[ \sigma =  { 7(1-\e^2)\over 4 ( \e + x^2+y^2 )^2  }
 + { 2\e \over \e + x^2+y^2 }.\]
That is, 
the Gauss curvature  is given by
\begin{eqnarray*}
{\bf K} & = &  3 {c_x u+ c_y v \over F} + \sigma \\
& = & -{ 3\sqrt{1-\e^2} (xu+yv)\over (\e +x^2+y^2)^2 F } +  { 7( 1-\e^2)\over 4 ( \e + x^2+y^2 )^2} +
{2\e \over \e +x^2+y^2}.
\end{eqnarray*}
For any Randers metric $F=\alpha+\beta$, the distortion is given by
\[ \tau = \ln \Big [ { F \over \alpha}\cdot {1 \over  1-\|\beta\|^2_{\alpha} } \Big ]^{n+1\over 2}.\]
A direct computation yields
\[ 1 - \|\beta\|^2_{\alpha} = { \e \Big ( 1+ \e (x^2+y^2)\Big ) \over \e + x^2+y^2 }.\]
Then the function $\nu: = \Big ( {\bf K} + {3c^2+\sigma\over 2} \Big )e^{ {2\tau\over n+1}} $ in (\ref{Kmu}) is given  by
\[ \nu = {3\over \e (\e +x^2+y^2) }.\]
That is, the Gauss curvature can also be given by
\begin{eqnarray*}
{\bf K} & = & - { 3c^2+\sigma\over 2} +\nu\;  e^{ -{2\tau\over n+1} }\\
& = & - { 5-\e^2 + 4\e (x^2+y^2)\over 2 (\e+x^2+y^2)^2 } 
+ { 3 \Big ( 1 +\e (x^2+y^2) \Big )\alpha 
\over (\e + x^2+y^2)^2 F } .
\end{eqnarray*}

}
\end{ex}

\section{Proof of Theorem \ref{ChenShen_thm1.1}}
\bigskip

Let $F = \alpha+\beta$ be a Randers metric on an $n$-dimensional manifold $M$,
where $\alpha = \sqrt{a_{ij}(x)y^iy^j}$ and $\beta =b_i(x)y^i$. 
Throughout this paper, we always assume that $F$ is  positive definite or  $\|\beta\|_{\alpha}(x) :
=\sqrt{a^{ij}(x) b_i(x) b_j(x)} < 1$ for any $x\in M$. 
Define $b_{i|j}$ by
\[ b_{i|j} \theta^j := db_i -b_j \theta_i^{\ j},\]
where 
$\theta^i :=dx^i$ and 
$\theta_i^{\ j} :=\gamma^i_{jk}dx^k$ denote the Levi-Civita  connection forms of $\alpha$.
Let
\[ r_{ij} :={1\over 2} \Big ( b_{i|j}+b_{j|i}\Big ) , \ \ \ \ \ 
s_{ij}:= {1\over 2} \Big (b_{i|j}-b_{j|i} \Big ),\]
\[s^i_{\ j}:= a^{ih}s_{hj}, \ \ \ \ \  s_j:=b_i s^i_{\ j}, \ \ \ \ \   e_{ij} := r_{ij}
+ b_i s_j + b_j s_i.\]
Let 
\[ \rho(x) := \ln\sqrt{1-\|\beta\|_{\alpha}^2(x) }\]
and $d\rho=\rho_i dx^i$. 
 According to \cite{Sh1}, the S-curvature of $F=\alpha+\beta$
is given by
\[
{\bf S} = (n+1) \Big \{ {e_{00} \over 2F}  - (s_0+\rho_0)   \Big \}, 
\]
where $e_{00}:=e_{ij}y^iy^j$ and $s_0:=s_iy^i$ and $\rho_0 := \rho_py^p$.
According to  Lemma 3.1 in \cite{ChSh}, ${\bf S}= (n+1)c(x) F$ is equivalent to  that 
\be
e_{ij}  = 2c(x) (a_{ij}-b_ib_j). \label{e_00}
\ee

Assume that $\alpha$ is of constant sectional curvature and $\beta$ is closed (hence $s_{ij}=0$ and $s_i=0$). 
Let 
\[ \Phi := b_{i|j} y^iy^j, \ \ \ \ \ \ \ \Psi:= b_{i|j|k} y^iy^jy^k.\]
By (8.56) in \cite{Sh1}, we have 
\be
{\bf K}F^2 =  \mu \alpha^2 + 3 \Big [ {\Phi\over 2F} \Big ]^2 - {\Psi \over 2F}.\label{Ric1}
\ee
Further we assume  that ${\bf S}= (n+1) c(x) F$. Since $s_{ij}=0$, $e_{ij}
=r_{ij} = b_{i|j}$ and 
(\ref{e_00}) simplifies to
\[ b_{i|j} = 2c (a_{ij}-b_ib_j).\]
We obtain 
\begin{eqnarray*}
\Phi & = & 2 c (\alpha^2-\beta^2)\\
\Psi & = & 2 c_{x^k} y^k (\alpha^2-\beta^2) - 8 c^2 (\alpha^2-\beta^2) \beta.
\end{eqnarray*}

\bigskip

Now we are ready to prove Theorem \ref{ChenShen_thm1.1}.
Let $F=\alpha+\beta$ be a Randers metric in Theorem \ref{ChenShen_thm1.1}.
Since $F$ is locally projectively flat, 
$\alpha$ is locally projectively flat and $\beta$ is closed \cite{MoSh}.
By the Beltrami theorem, we know that $\alpha$ is of constant sectional curvature $\mu$. Our main task is to determine $\beta$.

 By 
Theorem \ref{thmKS}, we know that the flag curvature is in the following form 
\be
{\bf K} = { 3 c_{x^k}(x) y^k \over F(x, y)} +\sigma(x), \label{KcF}
\ee
where $\sigma(x)$ is a scalar function on $M$. 
It follows from  (\ref{Ric1}) and (\ref{KcF})  that 
\[
 3 c_{x^k}y^k F +\sigma F^2 = {\bf K} F^2  = \mu \alpha^2 + 3 \Big [ {\Phi\over 2F} \Big ]^2 - {\Psi \over 2F} .
\]
Using the above formulas for $\Phi$ and $\Psi$, we obtain 
\[
2 \Big \{ 2 c_{x^k}y^k + (\sigma+c^2)\beta \Big \}
\alpha + \Big \{ 2c_{x^k}y^k + (\sigma+c^2) \beta \Big \} \beta + \Big \{ \sigma - 3 c^2-\mu \Big \} \alpha^2 =0.
\]
This gives
\begin{eqnarray}
&& 2 c_{x^k} y^k + (\sigma+c^2)\beta =0, \label{ccc1}\\
&& \sigma - 3 c^2 -\mu =0.\label{ccc2}
\end{eqnarray}
Plugging (\ref{ccc2}) into  (\ref{KcF}) and (\ref{ccc1})  yields
\be
{\bf K} = 3 \Big \{ {c_{x^k}(x)y^k\over F(x, y) } + c(x)^2 \Big \} +\mu.\label{KcF1}
\ee 
\be
2 c_{x^k} y^k + (\mu + 4 c^2 ) \beta =0. \label{ccc1*}
\ee
Now we are ready to determine $\beta$ and $c$. 

\bigskip
\noindent{\bf Case 1}: Suppose that $\mu + 4 c(x)^2 \equiv 0$. Then
$ c(x) = c$ is a constant. It follows from (\ref{KcF1}) that 
\[ {\bf K} =  3 c^2 +\mu = - c^2 .\]
Then Theorem \ref{ChenShen_thm1.1} (A) follows from the classification theorem for projectively flat Randers of constant curvature \cite{Sh2}.

\bigskip
\noindent{\bf Case 2}: Suppose that $\mu + 4 c(x)^2 \not=0$ on an open subset ${\cal U}\subset M$.  It follows from (\ref{ccc1*}) that 
\be
\beta = - { 2 c_{x^k} (x)y^k \over \mu + 4 c(x)^2 }. \label{betacy}
\ee 
Note that $\beta$ is exact. 
Let $c_i dx^i:= dc$ and $c_{i|j}  dx^j := dc_i - c_k \bar{\Gamma}^k_{ij} dx^j $  denote  the covariant derivative of $dc$ with respect to $\alpha$, were $\bar{\Gamma}^k_{ij}$ denote the Christoffel symbols of $\alpha$.  We have
\[ c_i = c_{x^i}(x), \ \ \ \ \ c_{i|j} = c_{x^ix^j}(x) -c_{x^k}(x)\bar{\Gamma}^k_{ij}(x) .\]
Similarly, we can define $b_{i|j}$ and $b_{i|j|k}$.   
Since $\beta$ is closed, $b_{i|j}= b_{j|i}$.
In this case, ${\bf S} = (n+1) c(x) F$ is equivalent to 
\be
b_{i|j} = 2 c (a_{ij} - b_i b_j). \label{babb}
\ee
From (\ref{betacy}), we have 
\be
b_i = - {2c_i\over \mu + 4 c^2}.\label{babb1}
\ee
Plugging (\ref{babb1}) into (\ref{babb}) yields
\be
c_{i|j} = - c (\mu + 4 c^2) a_{ij} + {12 c c_i c_j \over \mu + 4 c^2 }. \label{c_eq}
\ee

Next we are going to solve (\ref{c_eq}) for $c(x)$ in three cases when $\mu=-1, 0, 1$. 

\bigskip
\noindent{(B1)}: $\mu = -1$. We  assume that 
$\alpha=\alpha_{-1} = \sqrt{a_{ij}(x) y^iy^j}$ which is expressed in the form (\ref{alpha-1}). 
We have 
\[ a_{ij} = { \delta_{ij}\over 1-|x|^2} + {x^ix^j \over (1-|x|^2)^2 }.\]
The Christoffel symbols of $\alpha$ are given by
\[ \bar{\Gamma}^k_{ij} = {x^i \delta^k_j +  x^j \delta^k_i \over 1- |x|^2 }.\] 
Equation (\ref{c_eq}) becomes
\be
c_{x^ix^j}
- { x^i c_{x^j}+ x^j c_{x^i} \over 1-|x|^2}
= - c (-1 + 4 c^2) \Big \{ {\delta_{ij} \over 1-|x|^2}
+ {x^ix^j \over (1-|x|^2)^2 } \Big \}
+ { 12 c c_{x^i} c_{x^j} \over -1 + 4 c^2 }. \label{c_eq1}
\ee
Let 
\[   f : =   { 2c \sqrt{1-|x|^2}  \over \sqrt{ \mp (- 1 +4 c^2) }},\]
where the sign depends on the value of $c$ such that
$ \mp (-1+4c^2) >0$.
Equation (\ref{c_eq1}) simplifies to
\[ f_{x^i x^j} =0.\]
We obtain that  $f = \langle a, x \rangle  + \lambda $, where $\lambda \in \R$ and $a\in \R^n$. Then we obtain 
\be
c = { \lambda + \langle a, x\rangle \over 2 \sqrt{ ( \lambda + \langle a, x \rangle )^2 \pm (1-|x|^2 ) } }.\label{cc33}
\ee
Plugging (\ref{cc33}) into (\ref{betacy}) yields
\[
\beta = { (\lambda +\langle a, x \rangle)\langle x, y \rangle + (1-|x|^2) \langle a, y \rangle \over (1-|x|^2) \sqrt{ (\lambda + \langle a, x \rangle )^2 \pm (1-|x|^2) } }\]
and
\be
 F=  { \sqrt{|y|^2-(|x|^2|y|^2-\langle x, y \rangle^2)}\over 1-|x|^2} + { (\lambda +\langle a, x \rangle)\langle x, y \rangle + (1-|x|^2) \langle a, y \rangle \over (1-|x|^2) \sqrt{ (\lambda + \langle a, x \rangle )^2 \pm (1-|x|^2) } }. \label{F-1}
\ee

By a direct computation ,
\[1- \|\beta \|^2_{\alpha}
= { (1-|x|^2) \Big \{ \pm 1 - ( |a|^2 -\lambda^2 ) \Big \} \over
( \lambda +\langle a, x \rangle )^2 \pm (1-|x|^2)  }.\]
Clearly,   $F =\alpha+\beta$ is a Randers metric on an open subset of ${\rm B}^n$ if and only if $ |a|^2-\lambda^2< \pm 1$. In this case,
$( \lambda +\langle a, x \rangle )^2 \pm (1-|x|^2) >0$ for any $x\in {\rm B}^n$. Thus $F$ can be extended to the whole ${\rm B}^n$. By (\ref{KcF1}), (\ref{cc33}) and (\ref{F-1}), we obtain 
\[
{\bf K} = - {3\over 4} {\pm (1-|x|^2)\over ( \lambda +\langle a, x \rangle )^2 \pm (1-|x|^2)  }\cdot  { F(x, -y)\over F(x, y) } - {1\over 4}.
\]         

\bigskip
\noindent{(B2)} $\mu =0$. We assume that 
$\alpha=\alpha_0 = \sqrt{ \delta_{ij} y^iy^j}$ which  is expressed in the form (\ref{alpha-0}). Equation (\ref{c_eq}) becomes 
\be
c_{x^ix^j} = - 4 c^3 \delta_{ij} + { 3 c_{x^i} c_{x^j}\over c}.\label{c_eq2}
\ee 
Let ${\cal U}:=\{ x\in \R^n \ | \ c(x)\not=0\}$ and
let 
\[ f = {1\over c^2}.\]
Equation (\ref{c_eq2}) simplifies to 
\be
f_{x^ix^j} = 8 \delta_{ij}.\label{fij}
\ee
We obtain 
\[ f = 4 ( k  + 2\langle a, x \rangle + |x|^2 ),\]
where $k \in \R$ and $a\in \R^n$ such that $f(x) >0$ for $x\in {\cal U}$. Then 
$c=\pm 1/\sqrt{f}$ is given by
\be
c = { \pm 1\over 2 \sqrt{ k  + 2\langle a, x \rangle + | x|^2 } }.\label{cc37}
\ee
Plugging (\ref{cc37}) into (\ref{betacy}) yields
\[
\beta =\pm  { \langle a, y \rangle + \langle x, y \rangle \over \sqrt{ k  + 2 \langle a, x \rangle + |x|^2 } },
\]
and
\be
 F= |y| \pm  { \langle a, y \rangle + \langle x, y \rangle \over \sqrt{ k  + 2 \langle a, x \rangle + |x|^2 } }.\label{F_0}
\ee
Note that 
\[ 1 - \|\beta\|_{\alpha}^2 = { k  - |a|^2 \over k  + 2 \langle a, x \rangle + |x|^2 }.\]
Clearly, $F=\alpha+\beta$ is a Randers metric on an open subset 
of $\R^n$ if and only if  $ |a|^2 < k$. In this case,
\[  k + 2 \langle a, x \rangle + |x|^2 \geq  k - |a|^2 + ( |a| - |x|)^2  >0, \ \ \ \ \ \forall x \in \R^n.\] 
Thus $F$ can be extended to the whole $\R^n$. 
By (\ref{KcF1}), (\ref{cc37}) and (\ref{F_0}),  we obtain 
\[
 {\bf K} ={ 3 \over 4 ( k  + 2 \langle a, x \rangle + |x|^2)} 
\cdot {F(x, -y)\over F(x, y) }>0 .
\]

\bigskip
\noindent{(B3)} $\mu =+1$. We assume that 
$\alpha = \alpha_{+1} =\sqrt{a_{ij}(x)y^iy^j}$ which is expressed in the form 
(\ref{alpha+1}). 
We have 
\[ a_{ij} = {\delta_{ij}\over 1+|x|^2} - { x^i x^j \over (1+|x|^2)^2}.\]
The Christoffel symbols of $\alpha$ are given by
\[\bar{\Gamma}^k_{ij} = - {x^i \delta^k_j +  x^j \delta^k_i \over 1+|x|^2 }.     \]
Equation (\ref{c_eq}) becomes 
\be
c_{x^ix^j} + { x^i c_{x^j} + x^j c_{x^i} \over 1+|x|^2} 
= - c (1+4 c^2) \Big \{ {\delta_{ij}\over 1+|x|^2} - { x^i x^j \over (1+|x|^2)^2}\Big \} + { 12 c c_{x^i} c_{x^j}\over 1 + 4 c^2}.      \label{c_eq3}
\ee
Let 
\[ f := {2 c \sqrt{ 1+ |x|^2} \over \sqrt{1 + 4 c^2} }.\]
Equation (\ref{c_eq3}) simplifies to $f_{x^ix^j}=0$. We obtain that $f = \e + \langle a, x \rangle $. 
Then we obtain 
\[
c = { \e + \langle a,  x \rangle \over 2 \sqrt{ 1 + |x|^2 - (\e +\langle a, x \rangle )^2 } }.
\]
Thus 
\[
\beta =  { ( \e + \langle a, x \rangle ) \langle x, y \rangle  - (1+|x|^2) \langle a, y \rangle 
\over (1+|x|^2) \sqrt{1 + |x|^2 - (\e +\langle a, x \rangle )^2 }  }. 
\]
and 
\[ F =  { \sqrt{|y|^2+(|x|^2|y|^2-\langle x, y \rangle^2)}\over 1+|x|^2} + { (\e +\langle a, x \rangle)\langle x, y \rangle - (1+|x|^2) \langle a, y \rangle \over (1+|x|^2) \sqrt{  (1+|x|^2) -(\e + \langle a, x \rangle )^2 } }.\]

By a direct computation, 
\[ 1-\| \beta \|^2_{\alpha} 
= {  (1+|x|^2) \Big \{ 1 - \e^2 -|a|^2 \Big \} \over 
 1 + |x|^2 - (\e +\langle a, x \rangle )^2 }.
\]
Thus $F=\alpha+\beta$ is a  Randers metric 
on some open subset of $\R^n$ if and only if   $ \e^2 +|a|^2 < 1 $.
In this case, $1 + |x|^2 - (\e +\langle a, x \rangle )^2 >0$ for all $x\in \R^n$. Thus $F$ can extended to the whole $\R^n$. 
By (\ref{KcF1}), we obtain 
\[ {\bf K} =
{3(1+|x|^2 )\over 4 \{ 1+|x|^2 - (\e +\langle a, x \rangle)^2\} } \cdot { F(x, -y)\over F(x, y) } + {1\over 4} > {1\over 4} .\]

\bigskip
From Theorem \ref{ChenShen_thm1.1},  we obtain some interesting  projectively flat Randers metrics with isotropic S-curvature.

\begin{ex}{\rm 
Let 
\[
F_{-}(x, y): = { \sqrt{(1-|x|^2) |y|^2 + \langle x, y \rangle^2} \sqrt{  (1-|x|^2) +\lambda^2} +\lambda \langle x, y \rangle \over (1-|x|^2)\sqrt{  (1-|x|^2) +\lambda^2 } } ,
\ \ \ \ \ \ \ y\in T_x {\rm B}^n,
\]
where $\lambda \in \R^n$ is an arbitrary constant. The geodesics of $F_{-}$ are straight lines in ${\rm B}^n$.  Thus $F$ is of scalar curvature. One can easily verify that
$F_{-}$ is complete in the sense that every unit speed geodesic of $F_{-}$ is  defined on $(-\infty, \infty)$.   
Moreover $F_{-}$ has strictly negative flag curvature 
${\bf K} \leq -{1\over 4}$. }
\end{ex}

\bigskip
\begin{ex}{\rm 
Let  
\[
F_0(x, y) 
:= { |y| \sqrt{1+|x|^2} +\langle x, y \rangle \over \sqrt{ 1+|x|^2} }, 
\ \ \ \ \ \ \ y\in T_x\R^n. 
\]
The geodesics of  $F_0$ are straight lines in  $\R^n$. Thus $F$ is of scalar curvature.  One can easily verify that $F_0$ is positively complete in the sense that every unit speed geodesic of $F_0$ is defined on $(a, \infty)$ for some $a\in \R$.   
Moreover  $F_0$  has positive flag curvature ${\bf K} >0$. }
\end{ex}

\section{Proof of Theorem \ref{ChenShen_thm1.2}}

In Theorem \ref{ChenShen_thm1.2}, the manifold $M$ is compact. 
Assume that $\mu + 4c^2(x) \not=0$ on some open subset of $M$.

When $\mu \not=0$, 
let 
\[ f(x):= {2c(x) \over \sqrt{ \pm (\mu + 4 c(x)^2)} } ,\]
where the sign is chosen so that $\pm (\mu + 4c^2) >0$. 
We have 
\[
  f_{|i|j} = -\mu f a_{ij}.
\]
This gives
\be
\Delta f = -n \mu\;  f. \label{mun}
\ee

When $\mu=0$,  we take 
\[ f(x):= {1\over c(x)^2}.\]
We have 
\[
f_{|i|j} = 8 a_{ij}.
\]
This gives
\be
\Delta f = 8n. \label{8n}
\ee

\bigskip
\noindent{\bf Case 1}: $\mu = -1$. Suppose that $ 1 - 4 c(x)^2\not=0$ on $M$. 
Integrating (\ref{mun}) yields
\[ \int_M |\nabla f |^2 dV_{\alpha} 
= -\int_M f \Delta f dV_{\alpha} 
= - n \int_M f^2 \; dV_{\alpha} \]
Thus $ f =0$. This implies that $c=0$ and $F=\alpha$ is Riemannian.

Suppose that $1-4c(x_o)^2 \not=0$ at some point $x_o\in M$. 
Let $(\tilde{M}, \tilde{x}_o)$ be the universal cover of $(M, x_o)$.
We may assume that $\tilde{M}$ is isometric to
$({\rm B}^n, \alpha_{-1})$ with $\tilde{x}_o$ corresponding to the origin. The Randers metric $F$ is lifted to  a complete Randers metric $\tilde{F}$  on $\tilde{M}= {\rm B}^n$.
$\tilde{F}$  is given by (\ref{F-1}). Let
$\tilde{c}(\tilde{x})$ be the lift of $c(x)$, which is given by (\ref{cc33}). Thus 
$ 1-4 \tilde{c}(\tilde{x})^2\not=0$ for all $\tilde{x}\in {\rm B}^n$.  This implies that  $1-4c(x)^2 \not =0$ for all $x\in M$. By the above argument, we see that $c=0$. Hence $F=\alpha$ is Riemannian by (\ref{eq5}).

Suppose that $1 - 4 c(x)^2 \equiv 0$. Then  the  lift $\tilde{F}$ of $F$ to the universal cover $\tilde{M}={\rm B}^n$ is given by (\ref{F_a}), hence  it is incomplete. This is impossible because of the compactness of $M$. 
 We also see that 
$F$ has negative constant flag curvature and bounded Cartan torsion,  hence it  is Riemannian according to Akbar-Zadeh's theorem  \cite{AZ}\cite{Sh5}\cite{Sh6}. Then $c(x)=0$. 
This is a contradiction again.

\bigskip
\noindent{\bf Case 2}: $\mu = 0$. Suppose that $c(x_o) \not=0$.
Let $\tilde{M}$ denote the universal cover of $M$. 
We may assume that $\tilde{M}=\R^n$ with the origin corresponding to $x_o$.
The Randers metric $F$ lifted to $\tilde{M}=\R^n$ is given by (\ref{F_0}).
Thus $c(x) \not=0$ for all $x\in M$. 
Integrating (\ref{8n}) over $M$ yields
\[ 0 = \int_M \Delta f dV_{\alpha} = 8n {\rm Vol}(M, \alpha).\]
This is impossible.
Therefore  $c(x) \equiv 0$. In this case, $F$ is a locally projectively flat Randers metric with flag curvature ${\bf K}=0$, hence it is locally Minkowskian by \cite{Sh2}.

\bigskip
\noindent{\bf Case 3}: $\mu = 1$.  
Note that $1+ 4c(x)^2 \not=0$ on $M$. Let 
\[
 f(x):= {2c(x)\over \sqrt{1+4c(x)^2}}. \label{fc}
\]
It follows from (\ref{mun}) that 
\be
 f_{|i|j} = - f a_{ij}. \label{faij}
\ee
This gives
\[ \Delta f = - n f.\]
Thus $f$ is an eigenfunction of $(M, \alpha)$ with $\max_{x\in M} |f|(x) < 1$. 
We can express 
\be
 F(x,y) = \alpha(x,y) - {2c_{x^k}(x) y^k\over 1+ 4 c(x)^2} 
= \alpha(x,y) - { f_{x^k}(x) y^k \over \sqrt{ 1 - f(x)^2 } }. \label{F41}
\ee
\be
{\bf K}(x,y) 
= 3 \Big \{ {c_{x^k}(x)y^k \over F(x,y)} + c(x)^2 \Big \} + 1 
= {3\over 4 (1- f(x)^2)} { F(x, -y)\over F(x, y) } + {1\over 4}. \label{Kf}
\ee

Using (\ref{faij}), one can verify that 
\[ \delta :=\sqrt{ |\nabla f |_{\alpha}^2(x)+ f(x)^2 }\]
is a constant.
Since $F$ is positive definite,  $\delta < 1 $. 

Let 
\[ \lambda(x) := \sup_{y\in T_xM} { F(x, -y)\over F(x, y)} .\]
Using $|\nabla f|^2_{\alpha}(x) = \delta^2 - f(x)^2$, we obtain 
\[ \lambda (x) = { \sqrt{ 1- f(x)^2} + \sqrt{\delta^2 - f(x)^2} 
\over \sqrt{ 1-f(x)^2} - \sqrt{\delta^2 - f(x)^2} }.\]
Let $\lambda:=\max_{x\in M} \lambda(x)$. 
We have
\[ 1\leq \lambda(x) \leq \lambda = { 1+\delta \over 1-\delta}\]
and
\[ 1 -f(x)^2 = { (1-\delta^2) (\lambda(x) +1)^2 \over 4 \lambda(x)} .\]
Note that $\lambda(x)=\lambda$ if and only if $f(x)=0$. 
It follows from (\ref{Kf}) that 
\[
{2-\delta \over 2(1+\delta)}= {3+\lambda\over 4\lambda}  \leq {\bf K} \leq {3\lambda+1\over 4} = { 2 +\delta \over 2 (1-\delta)}.
\] 

Let 
\[ h(x):= \arctan \Big ( 2 c(x) \Big ).\]
The Randers metric $F(x, y)$ in (\ref{F41})  can be expressed by 
\[ F(x, y) = \alpha(x, y) - h_{x^k}(x) y^k.\]
Clearly  $F$ is pointwise projectively equivalent to $\alpha$, namely
the geodesics of $F$ are geodesics of $\alpha$ as point sets. Let $\sigma(t)$ be
a closed geodesic of $\alpha$. 
Observe that 
\[ F\Big (\sigma(t), \dot{\sigma}(t)\Big  )
=  \alpha \Big (\sigma(t), \dot{\sigma}(t)\Big  ) - {d\over dt} \Big [ h( \sigma(t)) \Big ] .\]
By the above equation we obtain 
\be
 {\rm Length}_F (\sigma) =\int F\Big (\sigma(t), \dot{\sigma}(t)\Big  )
  dt = \int   \alpha \Big (\sigma(t), \dot{\sigma}(t)\Big  )  dt = {\rm Length}_{\alpha} (\sigma).\label{sigma}
\ee
Assume that $M$ is simply connected. Then $(M, \alpha) ={\rm S}^n$.
Let $\sigma$ be an arbitrary great circle on ${\rm S}^n$. 
By (\ref{sigma}), 
\[ {\rm Length}_F (\sigma) = 2\pi.\]

\noindent
Xinyue Chen\\
Department of Mathematics, Chongqing Institute of Technology, Chongqing 400050,  P.R. China\\
chenxy58@163.net

\smallskip

\noindent
Xiaohuan Mo\\
LMAM, School of Mathematical Sciences, Beijing University, Beijing 100871, P.R. China\\
moxh@pku.edu.cn

\smallskip

\noindent
Zhongmin Shen\\
Department of Mathematical Sciences, Indiana University-Purdue University Indianapolis, 402 N. Blackford Street, Indianapolis, IN 46202-3216, USA.  \\
zshen@math.iupui.edu

\end{document}